\documentstyle[12pt]{article}
\newcommand{\sect}[1]{\section{#1}\setcounter{equation}{0}}

\font\mbn=msbm10 scaled \magstep1
\font\mbs=msbm7 scaled \magstep1
\font\mbss=msbm5 scaled \magstep1
\newfam\mbff
\textfont\mbff=\mbn
\scriptfont\mbff=\mbs
\scriptscriptfont\mbff=\mbss\def\mbf{\fam\mbff}
\def\Re{{\mbf R}}

\def\Z{{\mbf Z}}
\def\Co{{\mbf C}}

\def\Di{{\mbf D}}
\newtheorem{Th}{Theorem}[section]
\newtheorem{Lm}[Th]{Lemma}

\newtheorem{D}[Th]{Definition}
\newtheorem{Proposition}[Th]{Proposition}
\newtheorem{R}[Th]{Remark}
\author{Alexander Brudnyi\thanks{Research supported in part by NSERC.
}\\
Department of Mathematics\\
University of Toronto, Toronto\\
Canada}
\title{MATRIX-VALUED CORONA THEOREM FOR MULTIPLY CONNECTED DOMAINS} 
\date{January 17, 1999} 
\begin{document} 
\maketitle
\begin{abstract}
{Let $D\subset\Co$ be a bounded domain, whose boundary $B$ consists of 
$k$ simple closed continuous curves and $H^{\infty}(D)$ be the algebra
of bounded analytic functions on $D$. We prove the matrix-valued corona
theorem for matrices with entries in $H^{\infty}(D)$.
}
\end{abstract}
\sect{\hspace*{-1em}. Introduction.}
Let $D\subset\Co$ be a bounded domain, whose boundary $B$ consists of 
$k$ simple closed continuous curves.
Let $H^{\infty}(D)$ be the uniform algebra of bounded analytic 
functions on $D\subset\Co$ with the pointwise multiplication and with the 
norm
$$
||f||=\sup_{z\in D}|f(z)|.
$$
The maximal ideal space of $H^{\infty}(D)$ is defined by
$$
M(H^{\infty}(D))=
\{\phi\ : \phi\in Hom(H^{\infty},\Co),\ \phi\neq 0\}
$$
equipped with the weak $*$ topology induced by the dual space of $H^{\infty}$.
It is a compact Hausdorff space. A function $f\in H^{\infty}(D)$ 
can be thought of as a continuous function on $M(H^{\infty}(D))$ 
via the Gelfand transform 
$\hat f(\phi)=\phi(f)$\ $(\phi\in M(H^{\infty}(D)))$.
It is well known (see, e.g. [St]) that
if $f_{1},f_{2},...,f_{n}$ are functions in $H^{\infty}(D)$ such that
$$
|f_{1}(z)|+|f_{2}(z)|+...+|f_{n}(z)|\geq\delta>0\ \ \ \ \ {\rm for\ all}\ 
z\in D
$$
then there are $H^{\infty}(D)$ functions 
$g_{1},g_{2},...,g_{n}$ so that
$$
f_{1}(z)g_{1}(z)+f_{2}(z)g_{2}(z)+...+f_{n}(z)g_{n}(z)=1.
$$
This result is a generalization of the famous Carleson's corona theorem ([C])
and is equivalent to the statement that $D$ is dense in 
$M(H^{\infty}(D))$. In this paper we consider the following generalization
of Carleson's theorem.

Let $f=(f_{ij})$ be a $(k\times n)$-matrix, $k\leq n$, whose entries belong
to $H^{\infty}(D)$ and $\{F_{s}\}_{s\in S}$ be the
family of minors of $f$ of order $k$. Assume that
\begin{equation}\label{matr}
\sum_{s\in S}|F_{s}(z)|\geq\delta>0\ \ \ \ \ {\rm for\ all}\ z\in D .
\end{equation}
\begin{Th}\label{carl}
There exists a unimodular $(n\times n)$-matrix
$\widetilde f=(\widetilde f_{ij})$,  $\widetilde f_{ij}\in H^{\infty}$ 
for all $i,j$, so that $\widetilde f_{ij}=f_{ij}$ 
for $1\leq i\leq k$. 
\end{Th}
Our proof heavily relies upon the fact that the topological dimension of 
$M(H^{\infty}(D))$ equals 2, that is a consequence of the similar result
of Su'{a}rez ([S]) for $M(H^{\infty}(\Di))$, where $\Di\subset\Co$ is the unit 
disc.
\sect{\hspace*{-1em}. Proof of Theorem \ref{carl}.}
We begin with
\begin{Lm}\label{ru}
Let $B_{1},...,B_{k}$ be the components of the boundary of $D$, with
$B_{1}$ forming the outer boundary. Let $D_{1}$ be the interior of $B_{1}$,
and $D_{2},...,D_{k}$ the exteriors of $B_{2},...,B_{k}$, including the point
at infinity. For every $f\in H^{\infty}(D)$ there exists a decomposition
\begin{equation}\label{deco}
f(z)=f_{1}(z)+...+f_{k}(z)\ \ \ \ \ (z\in D),
\end{equation}
such that $f_{n}\in H^{\infty}(D_{n})$ for $n=1,...,k$.
\end{Lm}
For similar but more general statement see [R,\ Th. 3.3].
Fot the sake of completeness we present here a simple proof of the lemma.\\
{\bf Proof.}
By definition each $B_{i}$ can be approximated by a sequence
$\{C_{ip}\}_{p\geq 1}$ of simple smooth 
closed curves containing in $D$. Then $C_{1p},...,C_{kp}$ 
for sufficiently big $p$ are nonintersecting  smooth closed 
curves in $D$ which bound a domain $D_{p}\subset D$.
Put
$$
f_{np}(z):=\frac{1}{2\pi i}\int_{C_{np}}\frac{f(w)}{w-z}dw\ \ \ \ \
(n=1,...,k;\ z\in D_{p}).
$$
Then (\ref{deco}) holds in $D_{p}$. Moreover, we can choose $D_{p}$
so closed to $D$ that the following inequalities
$$
||f_{s}|_{C_{np}}||_{L^{\infty}(C_{np})}\leq A ||f||\ \ \ \ \
(1\leq s,n\leq k,\ s\neq n)
$$
hold with $A$ depending only on the distances 
between curves $B_{1},...,B_{k}$. Here $||\cdot ||$ stands for the $H^{\infty}$
norm of $f$. This follows directly from estimates of
Cauchy's kernels restricted to  $C_{np}$. 
The above inequality and decomposition (\ref{deco}) combined imply
$$
||f_{np}||_{L^{\infty}(C_{np})}\leq A'||f||
$$
for any $n=1,...,k$ with an absolute constant $A'$. Letting 
$p\to\infty$, one easily sees that there is 
$f_{n}=\lim_{p\to\infty}f_{np}$ such that $f_{n}\in H^{\infty}(D_{n})$,
$1\leq n\leq k$, and
$$
f(z)=f_{1}(z)+...+f_{k}(z)\ \ \ \ \ (z\in D)
$$
which completes the proof of the lemma.\ \ \ \ \ $\Box$

By Osgood-Caratheodory theorem for each $i=1,...,k$
there is a biholomorphic map $T_{i}:D_{i}\longrightarrow\Di$ 
that can be extended to continuous map 
$\overline{D_{i}}\longrightarrow\overline{\Di}$. Thus we obtain
isomorphism 
$T_{i}^{*}: H^{\infty}(\Di)\longrightarrow H^{\infty}(D_{i})$
which maps $H^{\infty}(\Di)\cap C(\overline{\Di})$ into
$H^{\infty}(D_{i})\cap C(\overline{D_{i}})$. Let $M(H^{\infty})$
be the maximal ideal space of $H^{\infty}(\Di)$ and
$M(H^{\infty}(D_{i}))$ be the maximal ideal space of 
$H^{\infty}(D_{i})$, $i=1,...,k$. Then $T_{i}$ can be extended
to homeomorphism $T_{i}':M(H^{\infty}(D_{i}))\longrightarrow
H^{\infty}(\Di)$.
\begin{Lm}\label{ext}
Let $f\in H^{\infty}(D)$. Then $f$ admits a continuous extension to
$\penalty-10000$
$(M(H^{\infty}(D_{i}))\setminus D_{i})\cup D$ for $i=1,...,k$.
\end{Lm}
{\bf Proof.}
First, notice that for any $i=1,...,k$ there is a natural 
continuous surjective mapping
$P_{i}:M(H^{\infty}(D_{i}))\setminus D_{i}\longrightarrow
B_{i}$ defined by the embedding homomorphism
$H^{\infty}(D_{i})\cap C(\overline{D_{i}})\longrightarrow
H^{\infty}(D_{i})$. Further,
according to Lemma \ref{ru}, $f=f_{1}+...+f_{k}$ with 
$f_{n}\in H^{\infty}(D_{n})$ for all $n$. Moreover,
$f_{s}|_{D_{i}}$ is continuous for $s\neq i$. Let 
$\{z_{\alpha}\}$ be a net in $D$ converging to a point
$\xi\in M(H^{\infty}(D_{i}))\setminus D_{i}$.
In particular, $\{z_{\alpha}\}$ converges to a point 
$P_{i}(\xi)\in B_{i}$ in the topology defined on $\Co$.
Since $f_{s}$, $s\neq i$, is continuous on $B_{i}$ 
we obtain (for such $s$)
$$
\lim_{\alpha}f_{s}(z_{\alpha})=f_{s}(P_{i}(\xi))\ .
$$
We stress that this limit does not depend on the choice of a net 
converging to $\xi$. 
Let us define
$$
f(\xi)=f_{i}(\xi)+\sum_{s\neq i}f_{s}(P_{i}(\xi))
$$
where by definition 
$$
f_{i}(\xi)=\lim_{\alpha}f_{i}(z_{\alpha}).
$$
It is easy to see that such extension of $f$ determines
a continuous function on 
$(M(H^{\infty}(D_{i}))\setminus D_{i})\cup D$. Note also
that the algebra of extended functions separates points of
$(M(H^{\infty}(D_{i}))\setminus D_{i})\cup D$.

The proof is complete.\ \ \ \ \ $\Box$

Let us consider now compact topological space
$$
K:=D\cup [\bigcup_{i=1}^{k}(M(H^{\infty}(D_{i}))\setminus D_{i}]
$$
with the topology induced from the corresponding maximal ideal
spaces. Then according to Lemma \ref{ext}, $H^{\infty}(D)$
admits a continuous extension to $K$ and the extended algebra
separates points of $K$. Since $D$ is dense in $K$ the corona
theorem for $H^{\infty}(D)$ obviously implies
that $K$ is homeomorphic to $M(H^{\infty}(D))$.
\begin{D}
For a normal space $X$ we say that $dim\ X\leq n$ if every finite
open covering of $X$ can be refined by an open covering whose order $\leq
n+1$. If $dim\ X\leq n$ and the statement $dim\ X\leq n-1$ is false,
we say that $dim\ X=n$. For the empty set we put $dim\ \emptyset=-1$.
\end{D}
\begin{Lm}\label{dimen}
$dim\ M(H^{\infty}(D))=2$.
\end{Lm}
{\bf Proof.}
We recall first
\begin{D}\label{mod}
Let $X$ be a normal space and $F\subset X$ be a closed subset. 
We say that 
$dim\ (X\ mod\ F)\leq n$ if $dim\ H\leq n$ for any closed set $H$ with 
$H\cap F=\emptyset$.
\end{D}
Our proof of the lemma rests upon the following result, 
see, e.g. [N,\ Ch.2 Th.9-11].
\begin{Proposition}\label{na}
Let $X$ be a normal space and $F$ be a closed subset of $X$. If 
$dim\ F\leq n$ and $dim\ (X\ mod\ F)\leq n$, then $dim\ X\leq n$.
\end{Proposition}
Set $F:=\cup_{i=1}^{k}(M(H^{\infty}(D_{i}))\setminus D_{i}$. 
Since by Definition \ref{mod}
$dim\ (M(H^{\infty}(D))\ mod\ F)=2$ and since  
$dim\ M(H^{\infty}(D_{i}))\setminus D_{i}\leq 2$ for $i=1,...,k$ due to
Su'{a}rez' theorem [S], application of
Proposition \ref{na} yields equality 
$dim\ M(H^{\infty}(D))=2$.\ \ \ \ \ $\Box$

We are now in position to prove Theorem \ref{carl}.

Let $f=(f_{ij})$ be a $(k\times n)$-matrix, $k\leq n$, with entries in
$H^{\infty}(D)$ and let the family of minors of $f$ of order $k$ satisfy
inequality (\ref{matr}).
According to the corona theorem we can extend
$f$ to $M(H^{\infty}(D))$ in such a way that the extended matrix satisfies 
(\ref{matr}) for all $\phi\in M(H^{\infty}(D))$. Now we will prove 
a general statement which implies immediately our result. 

Let $b$ be a $(k\times n)$-matrix whose entries belong to a commutative 
Banach algebra $A$ of complex-valued functions defined on maximal ideal 
space $M(A)$. Assume that $dim\ M(A)=2$. 
Assume also that condition (\ref{matr}) holds at each 
point of $M(A)$ for the family of minors of $b$ of order $k$ (this means
that these minors do not belong together to a maximal ideal). 

\begin{Proposition}
There is an invertible $(n\times n)$-matrix
with entries in $A$ which completes $b$.
\end{Proposition}
According to [L,\ Th.3], it suffices to complete $b$ in the
class of continuous matrix-functions on $M(A)$. Thus we have to
find an $(n\times n)$-matrix
$\widetilde b=(\widetilde b_{ij})$ with entries
$\widetilde b_{ij}$ in $C(M(A))$ with $det\ \widetilde b=1$ and 
$\widetilde b_{ij}=b_{ij}$ for $1\leq i\leq k$. 
Matrix $b$ determines trivial subbundle $\xi$ of complex rank 
$k$ in trivial vector bundle $\theta^{n}=M(A)\times\Co^{n}$. Let 
$\eta$ be an additional to $\xi$ subbundle of $\theta^{n}$. Then clearly
$\xi\oplus\eta$ is topologically trivial. We will prove that
$\eta$ is also topologically trivial. Then a trivialization
$s_{1},s_{2},...,s_{n-k}$ of $\eta$ will determine the required complement of
$b$.
\begin{Lm}\label{triv}
Let $X$ be a 2-dimensional Hausdorff compact and let $p:E\longrightarrow X$
be a locally trivial continuous vector bundle over $X$ of complex rank $m$.
If $m\geq 2$ then there is a continuous nowhere vanishing section of $E$.
\end{Lm}
{\bf Proof.}
Let $\gamma_{m}$ be canonical vector bundle over the Stiefel manifold
$V(m,k)$ of $m$-dimensional subspaces of $\Co^{k}$. (Recall that
$\gamma_{m}$ is defined as follows: the fibre over an $x\in V(m,k)$ is
the $m$-dimensional complex space representing $x$.)
According to the general theory of vector bundles [Hu,\ Ch.3, Th.5.5]
there is a continuous mapping $f$ of $X$ into 
$V(m,k)$ such that $f^{*}\gamma_{m}=E$.
Considering $V(m,k)$ as a compact submanifold of some $\Re^{s}$, $f$ can
be determined by a finite family of continuous on $X$ functions, 
$f=(f_{1},...,f_{s})$.
For $X'=f(X)$ let $U'$ be a compact  polyhedron containing $X'$ and such 
that $\gamma_{m}$ has a continuos extension to a vector bundle over $U'$. 
Let us consider the inverse limiting system determined by all possible
finite collections of continuous on $X$ functions such that the last $s$
functions of any such collection are coordinates of the mapping $f$
(in the same order). For any such family $f_{\alpha}:=
(f_{1},...,f_{\alpha-s},f)$ denote 
$X_{\alpha}=f_{\alpha}(X)\subset\Re^{\alpha}$. Further, let 
$p_{\alpha}^{\beta}:X_{\beta}\longrightarrow X_{\alpha}$,
$\beta\geq\alpha$, be the mapping induced by the natural projection
$\Re^{\beta}\longrightarrow\Re^{\alpha}$ onto the first 
$\alpha$ coordinates. Let $U_{\alpha}$ be a compact polyhedron containing
$X_{\alpha}$ defined in a small open neighbourhood of $X_{\alpha}$ so that
the inverse limit of the system $(U_{\alpha},p_{\alpha}^{\beta})$ coincides
with $X$ and $p_{0}^{\alpha}:\Re^{\alpha}\longrightarrow
\Re^{s}$ maps $U_{\alpha}$ to $U'$. Consider bundle
$E_{\alpha}:=(p_{0}^{\alpha})^{*}(\gamma_{m})$ over $U_{\alpha}$.  
Clearly, $(p_{\alpha}^{\beta})^{*}(E_{\alpha})=E_{\beta}$, $\beta\geq\alpha$,
and the pullback of each $E_{\alpha}$ to $X$ coincides with $E$. Now the
Euler class $e_{\alpha}$ of each $E_{\alpha}$ is an element of the
\v{C}ech cohomology group
$H^{2m}(U_{\alpha},\Z)$.
This class equals 0 if and only if $E_{\alpha}$ has a nowhere vanishing
section. By the fundamental property of characteristic classes
[Hi], $(p_{\alpha}^{\beta})^{*}(e_{\alpha})=e_{\beta}$, $\beta\geq\alpha$.
Further, it is well known [Br,\ Ch.II, Corol.14.6] that
$$
\lim_{\longrightarrow}(p_{\alpha}^{\beta})^{*}H^{k}(U_{\alpha},\Z)=
H^{k}(X,\Z),\ \ \ k\geq 0.
$$
Since $dim\ X=2$, Theorem 37-7 and Corollary 36-15 of [N] imply
$H^{k}(X,\Z)=0$ for $k>2$. But real rank of $E$ is $2m\geq 4$ so that
the Euler class $\displaystyle 
e=\lim_{\longrightarrow}(p_{\alpha}^{\beta})^{*}(e_{\alpha})$
of $E$ equals 0. From here and the above formula it follows
that there is some $\beta$ such that $e_{\beta}=0$. In particular,
$E_{\beta}$ has a continuous nowhere vanishing section $s$. Then its
pullback to $X$ determines the required section of $E$.\ \ \ \ \ $\Box$

Lemma \ref{triv} implies that $E$ is isomorphic to
direct sum $E_{m-1}\oplus E'$, where $E_{m-1}=X\times\Co^{m-1}$ and
$E'$ is a vector bundle over $X$ of complex rank 1. In particular, if
the first Chern class $c_{1}(E)\in H^{2}(X,\Z)$ of $E$ is 0, then 
$E'=X\times\Co$ and $E=X\times\Co^{m}$.
In our case $dim\ M(A)=2$, $\theta^{n}=\xi\oplus\eta$ and $\xi$ is
topologically trivial. Then for the first Chern class of $\theta^{n}$ we
have
$$
0=c_{1}(\theta_{n})=c_{1}(\xi)+c_{1}(\eta)=c_{1}(\eta).
$$
Therefore from Lemma \ref{triv} we obtain that $\eta$ is a 
topologically trivial bundle.

The proposition is proved.\ \ \ \ \ $\Box$\\
To complete the proof of the theorem it suffices to apply Lemma
\ref{dimen}.\ \ \ \ \ $\Box$
\begin{R} 
It is worth noting that similarly to [B] one can characterize topology of 
analytical part of $M(H^{\infty}(D))$.
\end{R}

\end{document}